\documentclass[preprint,12pt]{elsarticle}

\usepackage{amssymb}
\usepackage{algpseudocode}
\usepackage{algorithm}
\usepackage{bm}
\usepackage{listings}

\usepackage{lineno}
\usepackage{pgfplots}
\usepackage{pstricks,pst-node, pst-fun}
\usepackage{dsfont}
\usepackage{amssymb}
\usepackage[british,UKenglish]{babel}
\usepackage[utf8]{inputenc}
\usepackage[a4paper,left=3cm,right=2cm,top=2.5cm,bottom=2.5cm]{geometry}
\usepackage{amsmath}
\usepackage{relsize}
\usepackage{setspace}
\usepackage{subfig}
\DeclareMathAlphabet{\mathpzc}{OT1}{pzc}{m}{it}

\DeclareFixedFont{\ttb}{T1}{txtt}{bx}{n}{12} 
\DeclareFixedFont{\ttm}{T1}{txtt}{m}{n}{12}  

\usepackage{color}
\definecolor{deepblue}{rgb}{0,0,0.5}
\definecolor{deepred}{rgb}{0.6,0,0}
\definecolor{deepgreen}{rgb}{0,0.5,0}

\usepackage{listings}

\newcommand\pythonstyle{\lstset{
language=Python,
basicstyle=\ttm,
morekeywords={self,False,True},              
keywordstyle=\ttb\color{deepblue},
emph={MyClass,__init__},          
emphstyle=\ttb\color{deepred},    
stringstyle=\color{deepgreen},
frame=tb,                         
showstringspaces=false
}}

\lstnewenvironment{python}[1][]
{
\pythonstyle
\lstset{#1}
}
{}

\newcommand\pythoninline[1]{{\pythonstyle\lstinline!#1!}}

\journal{Journal Name}

\begin{document}

\begin{frontmatter}






\title{General form of the Gauss-Seidel equation to linearly approximate the Moore-Penrose pseudoinverse in random non-square systems and high order tensors}


\author{Luis Saucedo-Mora$^{a}$*, Luis Irastorza-Valera$^{a}$}

\address{$^a$ E.T.S. de Ingeniería Aeronáutica y del Espacio, Universidad Politécnica de Madrid, Pza. Cardenal Cisneros 3, 28040, Madrid, Spain}
\vspace{0.3cm}

\begin{abstract}
 
%
%
%
%
%
%

The Gauss-Seidel method has been used for more than 100 years as the standard method for the solution of linear systems of equations under certain restrictions. This method, as well as Cramer's and Jacobi's, is widely used in education and engineering, but there is a theoretical gap when we want to solve less restricted systems, or even non square or non-exact systems of equation. Here, the solution goes through the use of numerical systems, such as the minimization theories or the Moore-Penrose pseudoinverse. In this paper we fill this gap with a global analytical iterative formulation that is capable to reach the solutions obtained with the Moore-Penrose pseudoinverse and the minimization methodologies, but that analytically lies to the solutions of Gauss-Seidel, Jacobi, or Cramer when the system is simplified.

\end{abstract}


\end{frontmatter}



\section{Introduction} 

Optimization methods are widely used mathematical tools in Science and Engineering. Applications range from least squares regression---fitting data points \cite{takajo1988noniterative}---to training neural networks via backpropagation \cite{leung1991complex} and topology optimization in computational mechanics \cite{Bendsoe1989}. The choice of an optimization method depends on factors such as the problem's characteristics, convergence, and computational cost.

Optimization problems can be categorized as linear or non-linear, according to their objective function; and as constrained or unconstrained, depending on additional conditions. Convex problems are well-posed, i.e. their global optimum (solution) is easily reachable within a countable iteration number and distinguishable from local minima. However, most complex problems involving many variables and/or non-linearity are non-convex, whose many minima and saddle points may drive the optimization process away from their solution.

Analytical expressions yield exact solutions through a clearly defined model with known and measurable variables. If an analytical solution is available, it is preferred due to its deterministic nature. Alas, sometimes the analytical equations are either incomplete (lack of knowledge on the problem) or have too many variables (curse of dimensionality), on top of the bias induced by the user's choice of variables. Otherwise, numerical methods, which are iterative and require convergence criteria, are employed. Reducing the number of function evaluations in iterative methods minimizes computational cost, a key factor in optimization.

Numerical methods rely mostly on basic mathematical operations (multiplication, division, derivatives, etc.) performed upon data arranged in different arrays (tensors, vectors). In the matrix form, methods such as Jacobi, Gauss-Seidel and LU factorization are most prominent. If derivation is needed for vectors and fields (e.g. for differential equations in Boundary Value Problems), iterative methods such as Newton-Raphson and Gradient Descent are preferred. Other methods allow for discretization of variables to avoid differentiation: bisection, Runge-Kutta, Euler (forward/backward), Simpson, etc. 

If neither analytical nor numerical methods are up to the optimization task, other stochastic and/or statistic techniques can be used: the Montecarlo method, Machine Learning (particularly Neural Networks), etc. These options can prove really helpful when tackling highly dimensional and/or nonlinear problems whose analytical expressions are either computationally unmanageable or outright unavailable \cite{Jain2017,Sun2020,Dahrouj2021}. This sort of complex optimization problems can be simplified via Model Order Reduction (MOR) techniques \cite{Chinesta2011,Baur2014} such as Singular Value (SVD), Proper Orthogonal (POD) \cite{Lumley1967} or Proper Generalized Decomposition (PGD) \cite{Chinesta2014}; as well as Locally Linear Embeddings (LLE) \cite{Roweis2000}. 

However, these solutions are always approximated and the convergence time is often unknown \textit{a priori}. To solve optimization problems efficiently, one must account for the objective function's nature, the cost of evaluating it and computing its derivatives, the presence of constraints and many other factors. The modeler relies heavily on a large-enough database to train their ML surrogate to the desired accuracy level, along a careful choice of the optimizer (e.g. Stochastic Gradient Descent (SGD), adam\textcopyright \cite{Kingma2017}). This training process could be time-consuming and ineffective, hence the interest for numerical methods remains.

\section{Use of Linear Equation Systems}
Linear equations are common in the finite element method (FEM), turning differential equations governing physical problems into linear equations \cite{bathe2007finite}. In structural analysis, FEM reformulates the elasticity problem in terms of displacements at specific points, in which the stiffness matrix is typically symmetric, positive-definite, and often sparse when dealing with large-scale models. In linear systems (including square and overconstrained systems), the solution is often the minimum point of a quadratic function. The chosen optimization technique significantly impacts computational efficiency and accuracy.

\subsection{Effect of Noise in Machine Learning and Other Overconstrained Systems}
In overconstrained linear systems (more equations than unknowns), an exact solution may exist, but adding noise to the system introduces inconsistencies. Hence, finding a solution becomes an optimization problem, whose outcome is given by the chosen objective function and optimization algorithm. For a square (consistent) matrix system, all considered methods converge to a unique solution. 


\subsection{Analytical Solutions for Square Linear Systems}
Linear systems can be compatible (at least one solution) or incompatible (no solution) \cite{Burgstahler1983}. A system is compatible if the rank of its augmented matrix matches that of its coefficient matrix. To be uniquely solvable, the determinant of the coefficient matrix must be nonzero. Compatible systems with a unique solution are determinate, while those with infinite solutions are indeterminate. Undeterminate systems require additional boundary conditions for a specific solution. Several methods solve square linear systems analytically.

\subsection{Cramer’s Rule}
As a direct method, Cramer's rule provides an exact solution within a fixed number of steps, expressing each component as a ratio of determinants. It has been extended to non-square systems, rivaling Gaussian elimination under specific conditions \cite{Burgstahler1983}. However, its computational complexity grows factorially with system size, making it impractical for large systems where iterative or matrix decomposition methods are generally preferred.

\subsection{Iterative Methods}
In these techniques, after $k$ iterations (runs of a certain algorithm), an initial guess $\mathbf{x_0}$ converges (or not) to an \textbf{approximate} solution $\mathbf{x_k}$, within a given tolerance $\delta$ from the exact one, $\mathbf{x^*}$. See Equation \ref{eq:iterative_general}:

\begin{equation}
\bm{x}^{k+1} =\mathbf{M}\mathbf{x}^k+\mathbf{c}
\label{eq:iterative_general}
\end{equation}

If operator $\mathbf{M}$ and constant term $\mathbf{c}$ remain invariant, this is a fixed-point (stationary) method. If $\mathbf{M}$ is a matrix, its spectral radius (biggest absolute eigenvalues) must be less than 1 for convergence. Let the matrix equation system $\mathbf{A}x=b$ be considered from here onwards.

\subsubsection{Jacobi and Gauss-Seidel}
In Jacobi's stationary method, each row $i$ in the approximate solution $\mathbf{x^{k+1}}$ is a function of the elements in the same row in $\mathbf{A}$ and $\mathbf{b}$ and the previous step $\mathbf{x^k}$:
\begin{equation}
\bm{x_i}^{k+1} =\frac{1}{\bm{a_{ii}}}(\bm{b_i}-\sum_{j=1, i \neq j}^{n}\bm{a_{ij}}\bm{x_j}^k),\ \ \ \bm{i}=1,...,\bm{n}
\label{eq:Jacobi_terms}
\end{equation}

where $\bm{x_i}^1$ are the initial guesses ($\bm{x}^1$).




Gauss-Seidel method's improvement lies in computing new solution components $\mathbf{x_i}$ upon the already obtained ones, accelerating convergence:
\begin{equation}
\bm{x_i}^{k+1} =\frac{1}{\bm{a_{ii}}}(\bm{b_i}-\sum_{j=1, i \neq j}^{i-1}\bm{a_{ij}}\bm{x_j}^{k+1}-\sum_{j=i+1}^{n}\bm{a_{ij}}\bm{x_j}^k),\ \ \ \bm{i}=1,...,\bm{n}
\label{eq:Gauss_Seidel_terms}
\end{equation}




If $\mathbf{A}$'s diagonal is strictly dominant (for each row, $\left|\bm{a_{ii}}\right| \leq\sum_{j\neq i}^{n} \left|\bm{a_{ij}}\right|$), both Jacobi's and Gauss-Seidel's methods will converge to the solution. If $\mathbf{A}$ is symmetric ($\mathbf{A}$=$\mathbf{A^T}$) and positive definite ($\mathbf{z^T}\mathbf{A}\mathbf{z} > \mathbf{0}$ as long as $\mathbf{z} \neq \mathbf{0}$), Gauss-Seidel's method converges for whatever equation system. For Jacobi's method, an extra condition applies: 2$\mathbf{D}-\mathbf{A}$ must be positive definite as well.

\subsection{Inverse of higher order tensors and non-square nested matrices}

Higher order tensors are algebraic objects presenting order greater than 0 (scalars), 1 (vectors) and 2 (matrices), that is, $N$-dimensional arrays of data for $N>2$. They are used in Continuum Mechanics to study inhomogeneous, materially uniform, nonsimple gross bodies \cite{Morgan1976}. 

Unlike tensor addition, multiplication and inversion differ from matrices. For higher orders ($N>3$), they become recurrent: one must decompose them from order $N$ back to $N=2$ so they can be processed.
Inversion is a crucial step in many direct and iterative numerical methods. 

In this line, there are also inversion procedures for higher order isotropic (coordinates are invariant to rotation) tensors with minor symmetries (some local components are invariant to index permutation as well) \cite{Monchiet2010}.

\subsection{Principal optimization theories}

In this section we introduce some optimization theories that will be compared with the proposed algorithm, both linear and quadratic.

\subsubsection{Conjugate gradient for systems of linear equations}

The conjugate gradient is an iterative, first order gradient-based optimization method. Gradient descent yields $x_{k+1}$ by moving $x_k$ a certain $\alpha_k$ step along the gradient at $x_k$ until some convergence criteria are satisfied. For square $\mathbf{A}$ (positive-definite) matrices (even large sparse linear equations \cite{hestenes1952methods, fletcher1964function}), it performs well: 
the residuals are orthogonal, and search directions are conjugate (or $\mathbf{A}-$orthogonal i.e. $\bm{p}^i\cdot\mathbf{A}\bm{p}^j = 0, \forall i\neq j$), avoiding repetition of the same search directions (resulting Kyrlov subspace's property).

\subsubsection{Powell's method for systems of linear equations}

Powell's method is non-gradient-based. It starts from an initial point and a set of $n$ linearly independent search directions, and proofs that a quadratic is minimized after a number of iterations lesser or equal than $n$ \cite{powell1964efficient}. 
A step of magnitude $\alpha_i$ is performed on each direction to find the next point (1-D optimization along that direction $\bm{p}_i$). The search vector with the highest step is deleted from the set, and the linear combination $\sum_i^n\alpha_i\bm{p}_i$ is added as a new search direction. This will eventually turn all the directions in the set mutually conjugate, ensuring an exact minimum \cite{powell1964efficient}.

\subsubsection{BFGS for systems of linear equations}

The BFGS algorithm is a second-order iterative optimization method using local curvature i.e. the (approximated) Hessian to obtain the next search direction \cite{broyden1970convergence, fletcher1970new, goldfarb1970family, shanno1970conditioning}. It is analogous to Newton's method but computationally simpler ($\mathcal{O}(n^3)$ to $\mathcal{O}(n^2)$) since it does not require the Hessian matrix inversion every iteration, rather an approximate given by the Sherman-Morrinson formula \cite{sherman1950adjustment}:  $\mathbf{B}_{k+1}\left(\bm{x}_{k+1} - \bm{x}_{k}\right) = \nabla f(\bm{x}_{k+1}) - \nabla f(\bm{x}_k)$ \cite{broyden1970convergence, fletcher1964function, goldfarb1970family,shanno1970conditioning}.


\subsubsection{SLSQP for systems of linear equations}

The Sequential Least Squares Programming (SLSQP) iteratively addresses non-linear optimization problems admitting constraints within it, whose objective function needs to be differentiable twice \cite{kraft1988software}. For linear systems, the least squares problem is set as the objective function $f(\bm{x})$, and this latter requirement is fulfilled since the Hessian $\mathbf{B}(\bm{x})$ is always defined as $\mathbf{B}(\bm{x}) = \nabla^2f(\bm{x}) = 2\mathbf{A}^T\mathbf{A}$. Iteration $x_{k+1}$ is obtained from $x_{k}$ by moving along a search direction $\bm{p}_k$ a step of magnitude $\alpha_k$ i.e. $\bm{x}_{k+1} = \bm{x}_{k} + \alpha_k\bm{p}_k$. Its search direction minimizes the quadratic (local) approximation of the problem's Lagrange function. 

Note that the Newton's method is recovered when the problem is unconstrained, since the optimal search direction at every step $k$ is $\bm{p}_k = -\mathbf{B}_k^{-1}\nabla f(\bm{x}_k)$. Also note that the optimal step size is $\alpha_k = 1$ in that case \cite{fliege2009newton}.

\subsection{Moore-Penrose pseudoinverse}

The unique Moore-Penrose (pseudo)inverse of a matrix $\mathbf{A}$ is $\mathbf{A}^{\dagger}$ satisfying
\begin{equation}
\mathbf{A}\mathbf{A}^{\dagger}\mathbf{A}=\mathbf{A};     \mathbf{A}^{\dagger}\mathbf{A}\mathbf{A}^{\dagger}=\mathbf{A}^{\dagger};    (\mathbf{A}\mathbf{A}^{\dagger})^T=\mathbf{A}\mathbf{A}^{\dagger};  (\mathbf{A}^{\dagger}\mathbf{A})^T=\mathbf{A}^{\dagger}\mathbf{A}
\end{equation}

For square, full rank matrices, $\mathbf{A}^{\dagger}=\mathbf{A}^{-1}$. It can also be obtained via Singular Value Decomposition (SVD) - or EigenValue Decomposition (EVD) if $\mathbf{A}$ is square and symmetric:
\begin{equation}
\mathbf{A}^{\dagger}=\mathbf{U}\mathbf{\Sigma}^{\dagger}\mathbf{V}^T
\label{eq:Moore_Penrose_SVD}
\end{equation}



\section{Theoretical model}

\subsection{The kernel of a function as the subspace of maxima in an exponential function}

The kernel $\mathpzc{K}$ of a function $\mathpzc{g}(\textbf{x})$ is a subset of points of a space of variables $\mathds{R}^n$ where the image of the function is null. It means that for some functions this region is of high interest because it is the solution of a certain calculation. Also, when the function is defined positive, if it exists, the kernel is the minimization of the function. And, it is not necessarily unique, since some equations may have several zeros.\\

One of the main problems in order to find the points of the kernel is the non-smoothness of the base function $\mathpzc{g}$, but if instead of solving $\mathpzc{g}(\textbf{x})=0$ we solve $\mathpzc{f}(\textbf{x})=\exp\left[-\mathpzc{g}(\textbf{x}) \right]=1$, using the exponential function as a base, the smoothness is ensured for any function $\mathpzc{g}(\textbf{x})$, and the solutions in $\mathds{R}^n$ for both equations are the same. We will call this the exponential kernel of $\mathpzc{g}(\textbf{x})$.\\

So here we ensure to have a smooth function (derivable infinite times) where the set of maxima is the kernel set of the function. Also, if we define positive the function $\mathpzc{f}(\textbf{x})=\mathpzc{g}(\textbf{x})^2$, and thus the equivalent exponential is $\exp\left[-\mathpzc{g}^2 \right]=1$, the problem now is to find the maximum of an exponential function, i.e. where the derivative is null.\\

\begin{figure}[H]
\centering
\includegraphics[width=0.8\textwidth]{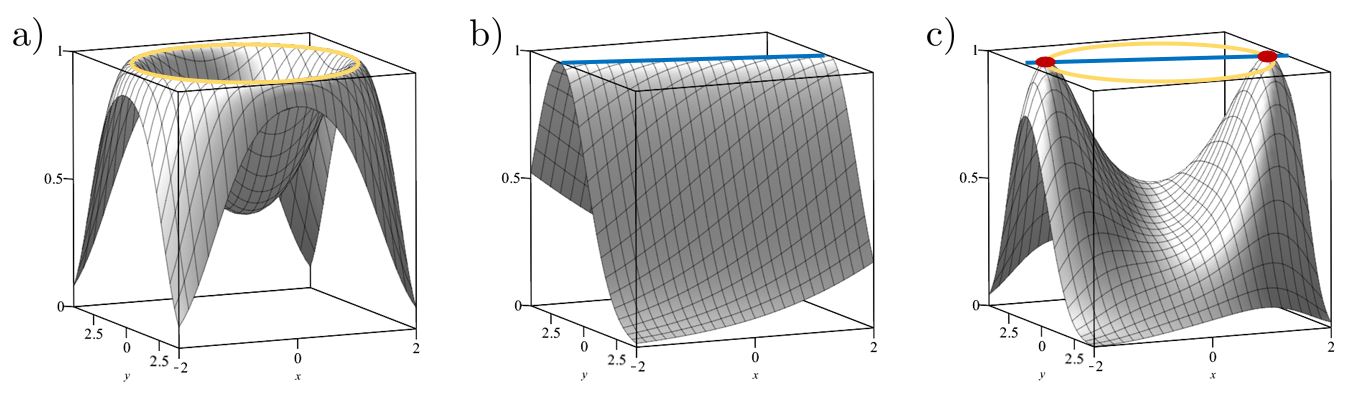}
\caption[tab]{The kernel of a function as the exponential subspace of maxima. For a circle in a), a line in b) and the combination of them in c).}
\label{fig:kernel}
\end{figure}

This is implemented in the example seen in Figure \ref{fig:kernel}, where the exponential kernel is used to solve the crossing points between a circle and a line. Let´s define the points of the circle through the expression $\mathpzc{c}(\textbf{x})=0$, and the points of the line through the equation $\mathpzc{l}(\textbf{x})=0$. So the crossing points between the circle and the line are the points that belongs simultaneously to the kernels of $\mathpzc{c}(\textbf{x})$ and $\mathpzc{l}(\textbf{x})$. With this, Figure \ref{fig:kernel}a shows the equation $\exp\left[-\mathpzc{c}(\textbf{x})^2 \right]$, and \ref{fig:kernel}b shows the equation $\exp\left[-\mathpzc{l}(\textbf{x})^2 \right]$. Then, the points that belongs to both kernels at the same time can be adressed through the multiplication of both exponential functions as $\exp\left[-\mathpzc{c}(\textbf{x})^2 \right]\,\exp\left[-\mathpzc{l}(\textbf{x})^2 \right]=\exp\left[-(\mathpzc{c}(\textbf{x})^2+\mathpzc{l}(\textbf{x})^2 ) \right]$, shown in Figure \ref{fig:kernel}c. \\

Finally, it is simple to point out that the set of solutions of $\exp\left[-(\mathpzc{c}(\textbf{x})^2+\mathpzc{l}(\textbf{x})^2 ) \right]=1$ is the same as the set of solutions of $\mathpzc{c}(\textbf{x})^2+\mathpzc{l}(\textbf{x})^2 =0$, which determines the crossing points between the line and the circle. This equivalence has the potentiality to be immersed in a smooth equation through the use of the exponential function. This feature will be used for the solution of systems of linear equations through the linear approximation of the gradient of the system of equations.\\

In general terms, a set of equations $\mathpzc{h}_i(\textbf{x})$ can be concatenated to find a subspace which simultaneously belongs to the kernel of all the functions:\\

\begin{equation}
\mathpzc{g}(\textbf{x})=\prod_{\hat{\textbf{i}}}\exp\left[-\frac{\mathpzc{h}_{\hat{\textbf{i}}}(\textbf{x})^2}{\gamma} \right]=\exp\left[-\sum_{\hat{\textbf{i}}} \frac{\mathpzc{h}_{\hat{\textbf{i}}}(\textbf{x})^2}{\gamma} \right] 
\end{equation}

And if the functions $\mathpzc{h}_i(\textbf{x})$ are systems of linear equations, it will have the form:\\

\begin{equation}
\mathpzc{g}(\textbf{x})=\prod_{\hat{\textbf{i}}} \exp \left[-\frac{(\mathcal{A}_{\hat{\textbf{i}}j}\textbf{x}_{j}-\textbf{b}_{\hat{\textbf{i}}})^2}{\gamma} \right]= \exp \left[-\sum_{\hat{\textbf{i}}} \frac{(\mathcal{A}_{\hat{\textbf{i}}j}\textbf{x}_{j}-\textbf{b}_{\hat{\textbf{i}}})^2}{\gamma} \right]  
\label{ecu:g1}
\end{equation}

Where $\hat{\textbf{i}}$ defines the equation evaluated from the system, and the parameter $\gamma$ is added for calculation purposes, and does not have any influence on the elements of the solution subset, only on the gradient of the exponential function. It will be used in the next section.\\

\subsection{Application to systems of linear equations}

The systems that will be analysed here come in the form:\\

\begin{equation}
\mathcal{A} \cdot \textbf{x}= \textbf{b}
\label{ecu:general}
\end{equation}

Where, according to equation \ref{ecu:g1}, $\mathpzc{g}(\textbf{x})$ is the sum of all the rows in equation \ref{ecu:general}, and equivalently we can operate to define: \\

\begin{equation}
\mathpzc{g}(\textbf{x})=\prod_{\hat{\textbf{i}}} \exp \left[-\frac{(\textbf{b}_{\hat{\textbf{i}}}-\mathcal{A}_{\hat{\textbf{i}}j}\textbf{x}_{j})^2}{\gamma} \right]= \exp \left[- \frac{(\textbf{b}_{i}-\mathcal{A}_{ij}\textbf{x}_{j})^T(\textbf{b}_{i}-\mathcal{A}_{ij}\textbf{x}_{j})}{\gamma} \right]
\label{ecu:g2}  
\end{equation}

Which is a smooth function. We derive this function to define the maximum of $\mathpzc{g}(\textbf{x})$, since the solution will be the point that makes the derivative null. The derivative with respect to a coordinate $\textbf{x}_{\hat{\textbf{j}}}$ of $\mathpzc{g}(\textbf{x})$ is:\\

\begin{equation}
\mathpzc{f}_{\hat{\textbf{j}}}(\textbf{x})=\frac{\partial \mathpzc{g}(\textbf{x})}{\partial x_{\hat{j}}}=-\frac{2 \mathcal{A}_{i\hat{\textbf{j}}}^T \left(\textbf{b}_{i}-\mathcal{A}_{ik}\textbf{x}_{k}\right)}{\gamma} \mathpzc{g}(\textbf{x}) =0
\label{deriv}
\end{equation}

It needs to be remarked that the derivative of equation \ref{deriv} is with respect to the value of the coordinate $x_{\hat{j}}$, not the axis $\hat{j}$ if we are thinking about a multidimensional space.\\

In equation \ref{deriv}, if $\gamma \rightarrow \infty$ the modulus of the derivative will be smaller as we approach to the solution point, then in this case $\mathpzc{g}(\textbf{x})\approx 1$, and:\\

\begin{equation}
\mathpzc{f}_{\hat{\textbf{j}}}(\textbf{x})=\mathcal{A}_{i\hat{\textbf{j}}}^T \left(\textbf{b}_{i}-\mathcal{A}_{ik}\textbf{x}_{k}\right) =-0\, \frac{\gamma}{2\mathpzc{g}(\textbf{x})}=0; \; \; \forall \textbf{x}\in \mathds{J}^n \in \mathds{R}^n
\label{deriv2}
\end{equation}

Where $\mathds{J}^n$ is the subspace of dimension $n$ composed by the solutions of the system.\\

The next step is to extract the value of $x_{\hat{j}}$ that fulfils equation \ref{deriv2}. For this we compute:\\   

\begin{equation}
\mathpzc{f}_{\hat{\textbf{j}}}(\textbf{x})=\mathcal{A}_{i\hat{\textbf{j}}}^T \left(\textbf{b}_{i}-\mathcal{A}_{ik}\textbf{x}_{k}(1-\delta_{\hat{\textbf{j}}k})-\mathcal{A}_{i\hat{\textbf{j}}}\,x_{\hat{j}}\right)=0 
\label{deriv}
\end{equation}

And,\\

\begin{equation}
\mathcal{A}_{i\hat{\textbf{j}}}^T \left(\textbf{b}_{i}-\mathcal{A}_{ik}\textbf{x}_{k}(1-\delta_{\hat{\textbf{j}}k})\right)=\mathcal{A}_{i\hat{\textbf{j}}}^T \mathcal{A}_{i\hat{\textbf{j}}} \,x_{\hat{j}}
\label{deriv}
\end{equation}

So,\\

\begin{equation}
x_{\hat{j}}=\frac{\mathcal{A}_{i\hat{\textbf{j}}}^T \left(\textbf{b}_{i}-\mathcal{A}_{ik}\textbf{x}_{k}(1-\delta_{\hat{\textbf{j}}k})\right)}{\mathcal{A}_{i\hat{\textbf{j}}}^T \mathcal{A}_{i\hat{\textbf{j}}}}; \;\; \forall x_{\hat{j}} \in \textbf{x}
\label{ecuxsyslin}
\end{equation}

In equation \ref{ecuxsyslin} the value of one of the coordinates of the solution is placed as a function of the rest of the coordinates of $\textbf{x}$. This is the general solution proposed in this paper for systems of linear equations, and in the next sections this will be analysed and extended to the solution of high order tensor and non-square nested matrices under dot product and double dot product.\\

\subsubsection{Iterative formulation}

The solution shown in equation \ref{ecuxsyslin} can be easily solved for a small system of equations, but has the drawback that for a large systems of equations, it is not possible to solve the system by extracting each variable independently. For this reason, we explore in this section the feasibility of the iterative resolution of the system.\\

For this, from Equation \ref{ecuxsyslin} we can operate with the $\delta_{\hat{\textbf{j}}k}$ term to obtain:\\

\begin{equation}
\mathcal{A}_{ik}\textbf{x}_{k}(1-\delta_{\hat{\textbf{j}}k})=\mathcal{A}_{ik}\textbf{x}_{k}-\mathcal{A}_{ik}\textbf{x}_{k}\delta_{\hat{\textbf{j}}k}=\mathcal{A}_{ik}\textbf{x}_{k}-\mathcal{A}_{i\hat{\textbf{j}}}\, x_{\hat{j}}
\label{xiter_desp}
\end{equation}

It implies that when the result of Equation \ref{xiter_desp} is substituted into Equation \ref{ecuxsyslin}, we obtain two $x_{\hat{j}}$ terms that will be used for the iterative formulation. Here we are going to differentiate between the three terms related to $\textbf{x}$ in the Equation \ref{ecuxsysliniter}, which conducts to the final formulation as:\\

\begin{equation}
\widehat{x}_{\hat{j}}=\frac{\mathcal{A}_{i\hat{\textbf{j}}}^T \left(\textbf{b}_{i}-\mathcal{A}_{ik}\breve{\textbf{x}}_{k}+\mathcal{A}_{i\hat{\textbf{j}}}\tilde{x}_{\hat{j}})\right)}{\mathcal{A}_{i\hat{\textbf{j}}}^T \mathcal{A}_{i\hat{\textbf{j}}}}; \;\; \forall \widehat{x}_{\hat{j}} \in \textbf{x}
\label{ecuxsysliniter}
\end{equation}

So, $\widehat{x}_{\hat{j}}$ is the term $\hat{j}$ at the new iteration $t+1$, $\tilde{x}_{\hat{j}}$ is the value at $\hat{j}$ of the previous iteration $t$, and the vector  $\breve{\textbf{x}}_{k}$ has all the values of the vector $\textbf{x}$. In the next section we will study the values of $\breve{\textbf{x}}_{k}$ regarding to $t+1$ and $t$.\\

\subsubsection{General convergence and error propagation}

We start from the solution of equation \ref{ecuxsyslin}, where we define $\mu_i^t$ as the error of the system when the vector $\widehat{x}$ is evaluated at iteration $t$, and it is not necessarily the system's solution. When the solution is reached, every $\mu_i^t$ at iteration $t_{sol}$ when the system has converged is equal to 0. Also, the system does not necessarily have a solution, due to possible noise in a non-square system. We denote this error as $\lambda_i$, which is the residual error of the system when the minimum solution of the system is reached. It is a vector of the same size of $b$, with the residual of every equation at its minimum, so it is constant over time.\\

For a certain iteration step $t$:\\

\begin{equation}
\boldsymbol{\mu}_i^t+\boldsymbol{\lambda_i}=\textbf{b}_{i}-\mathcal{A}_{ik}\textbf{x}_{k}^t(1-\delta_{\hat{\textbf{j}}k})-\mathcal{A}_{i\hat{\textbf{j}}}\,x_{\hat{j}}^t
\label{gamma1}
\end{equation}

So:\\

\begin{equation}
\mathcal{A}_{ik}\textbf{x}_{k}^t(1-\delta_{\hat{\textbf{j}}k})=\textbf{b}_{i}-\mathcal{A}_{i\hat{\textbf{j}}}\,x_{\hat{j}}^t-\boldsymbol{\mu}_i^t-\boldsymbol{\lambda_i}
\label{gamma2}
\end{equation}

Where $\boldsymbol{\lambda_i}$ is constant for all the calculation, since it is the error between the optimum solution and the values that fulfil equation $\mathcal{A}\textbf{x}=\textbf{b}$, so it is a corrector of the vector $\textbf{b}$ to make the system exact. The vector $\boldsymbol{\lambda_i}$ is null in a square system, but not necessarily in a non-square system with noise.\\ 

And inserting the error $\boldsymbol{\lambda}_i$ in the equation \ref{ecuxsyslin}:\\

\begin{equation}
x_{\hat{j}}=\frac{\mathcal{A}_{i\hat{\textbf{j}}}^T \left(\textbf{b}_{i}+\boldsymbol{\lambda}_i-\mathcal{A}_{ik}\textbf{x}_{k}(1-\delta_{\hat{\textbf{j}}k})\right)}{\mathcal{A}_{i\hat{\textbf{j}}}^T \mathcal{A}_{i\hat{\textbf{j}}}}; \;\; \forall x_{\hat{j}} \in \textbf{x}
\label{ecuxsyslinnu}
\end{equation}

And now substituting here equation \ref{gamma2} we have:\\

\begin{equation}
x_{\hat{j}}^{t+1}=\frac{\mathcal{A}_{i\hat{\textbf{j}}}^T \left(\mathcal{A}_{i\hat{\textbf{j}}}\,x_{\hat{j}}^t+\boldsymbol{\mu}_i^t     \right)}{\mathcal{A}_{i\hat{\textbf{j}}}^T \mathcal{A}_{i\hat{\textbf{j}}}}=x_{\hat{j}}^t  +\frac{\mathcal{A}_{i\hat{\textbf{j}}}^T \boldsymbol{\mu}_i^t }{\mathcal{A}_{i\hat{\textbf{j}}}^T \mathcal{A}_{i\hat{\textbf{j}}}}
\label{ecuxsyslin_gamma}
\end{equation}

Where the first term indicates the position of the point and the second the error that corrects the new iteration. Operating with 2 consecutive iterations starting from $t=t_0$:\\

\begin{equation}
x_{\hat{j}}^{t_0+2}=x_{\hat{j}}^{t_0}  +\frac{\mathcal{A}_{i\hat{\textbf{j}}}^T \boldsymbol{\mu}_i^{t_0} }{\mathcal{A}_{i\hat{\textbf{j}}}^T \mathcal{A}_{i\hat{\textbf{j}}}} +\frac{\mathcal{A}_{i\hat{\textbf{j}}}^T \boldsymbol{\mu}_i^{t_0+1} }{\mathcal{A}_{i\hat{\textbf{j}}}^T \mathcal{A}_{i\hat{\textbf{j}}}}
\label{ecuxsyslin_gamma2}
\end{equation}

And in a general form:\\

\begin{equation}
x_{\hat{j}}^{t}=x_{\hat{j}}^{t_0}  +\sum_{\xi=t_0}^{t} \frac{\mathcal{A}_{i\hat{\textbf{j}}}^T \boldsymbol{\mu}_i^{\xi} }{\mathcal{A}_{i\hat{\textbf{j}}}^T \mathcal{A}_{i\hat{\textbf{j}}}} 
\label{ecuxsyslin_gamma3}
\end{equation}

From here we can define:\\

\begin{equation}
\zeta_i^t=\sum_{\xi=t_0}^{t} \frac{\boldsymbol{\mu}_i^{\xi} }{\mathcal{A}_{i\hat{\textbf{j}}}^T \mathcal{A}_{i\hat{\textbf{j}}}} 
\label{ecuxsyslin_gamma5}
\end{equation}

So, the accumulative $\Delta x_{\hat{j}}$ from time $t_0$ to a general time $t=t_0+k$ will be:\\

\begin{equation}
\Delta x_{\hat{j}}^{t+1}=\sum_{t=t_0}^{t_0+k} \frac{\mathcal{A}_{i\hat{\textbf{j}}}^T \boldsymbol{\mu}_i^{t} }{\mathcal{A}_{i\hat{\textbf{j}}}^T \mathcal{A}_{i\hat{\textbf{j}}}} =\mathcal{A}_{i\hat{\textbf{j}}}^T \zeta_i^t
\label{ecuxsyslin_gamma6}
\end{equation}

Now, we can operate over $\Delta x_{\hat{j}}^{t}$ through the update of the error $\boldsymbol{\zeta}^t$ in different ways to prevent error propagation with an average of the values that potentially can lead towards a divergence of the method. For a simple example, we are going to denote $\alpha_{ij}$ as a component of the matrix $\mathcal{A}$, and  with this, in a system of 5 unknowns, the component 3 which is in the middle will be defined in the general form as:\\

\begin{equation}
\Delta x_3^{t+1}=\alpha_{31} \zeta_1^t+\alpha_{32} \zeta_2^t+\alpha_{33} \zeta_3^t+\alpha_{34} \zeta_4^t-\alpha_{35} \zeta_5^t
\label{ecu1}
\end{equation}

Which is a direct multiplication of the error vector by the components of the column of $\mathcal{A}$. Now we propose different scenarios for this example: \\

\textbf{1}.- The matrix $\mathcal{A}$ is symmetric we can write it as:

\begin{equation}
\Delta x_3^{t+1}=\alpha_{31} \left( \zeta_1^t+\zeta_5^t\right)+\alpha_{32} \left( \zeta_2^t+\zeta_4^t\right)+\alpha_{33}\zeta_3^t
\label{ecu2}
\end{equation}

\textbf{2}.- With a dynamic updating of the error we have:

\begin{equation}
\Delta x_3^{t+1}=\alpha_{31} \zeta_1^{t+1}+\alpha_{32} \zeta_2^{t+1}+\alpha_{33} \zeta_3^t+\alpha_{34} \zeta_4^t+\alpha_{35} \zeta_5^t
\label{ecu3}
\end{equation}

\textbf{3}.- With a $\beta$ coefficient to average the value of two consecutive iterations we define:

\begin{equation}
\eta_i^t=\beta\,\zeta_i^{t+1}\,+(1-\beta)\,\zeta_i^{t}
\end{equation}

Where without dynamic updating is:

\begin{equation}
\Delta x_3^{t+1}=\alpha_{31} \eta_1^t+\alpha_{32} \eta_2^t+\alpha_{33} \eta_3^t+\alpha_{34} \eta_4^t+\alpha_{35} \eta_5^t
\label{ecu4}
\end{equation}

\textbf{4}.- And with dynamic updating:

\begin{equation}
\Delta x_3^{t+1}=\alpha_{31} \eta_1^{t+1}+\alpha_{32} \eta_2^{t+1}+\alpha_{33} \eta_3^t+\alpha_{34} \eta_4^t+\alpha_{35} \eta_5^t
\label{ecu5}
\end{equation}

The equations \ref{ecu1}, \ref{ecu3}, \ref{ecu4}, \ref{ecu5} are different ways to solve the problem, with different treatments of error propagation. Each approach can handle the influence of the error of a certain component in the updated values of the rest of the variables. So, equation \ref{ecu1} has a free spread of the error, mitigated by the case of equation \ref{ecu2} when the matrix is symmetric. Equation \ref{ecu3} is an approach where the influence of the error is distributed through half of the variables in two consecutive iterations. And equations \ref{ecu4} and \ref{ecu5} incorporate the $\beta$ coefficient to gradually update the new value of the variable calculated.\\

Considering the more general approach, shown in equation \ref{ecu5}, equation \ref{ecuxsysliniter} can be reformulated in the form:\\

\begin{equation}
\widehat{\textbf{x}}_{\hat{\textbf{j}}}=(1-\beta)\tilde{\textbf{x}}_{\hat{\textbf{j}}}+\beta\frac{\mathcal{A}_{i\hat{\textbf{j}}}^T \left(\textbf{b}_{i}-\mathcal{A}_{ik}\breve{\textbf{x}}_{k}+\mathcal{A}_{i\hat{\textbf{j}}}\tilde{\textbf{x}}_{\hat{\textbf{j}}})\right)}{\mathcal{A}_{i\hat{\textbf{j}}}^T \mathcal{A}_{i\hat{\textbf{j}}}}; \;\; \forall \widehat{\textbf{x}}_{\hat{\textbf{j}}} \in \textbf{x}
\label{ecuxbeta}
\end{equation}

\begin{figure}[H]
\centering
\includegraphics[width=0.8\textwidth]{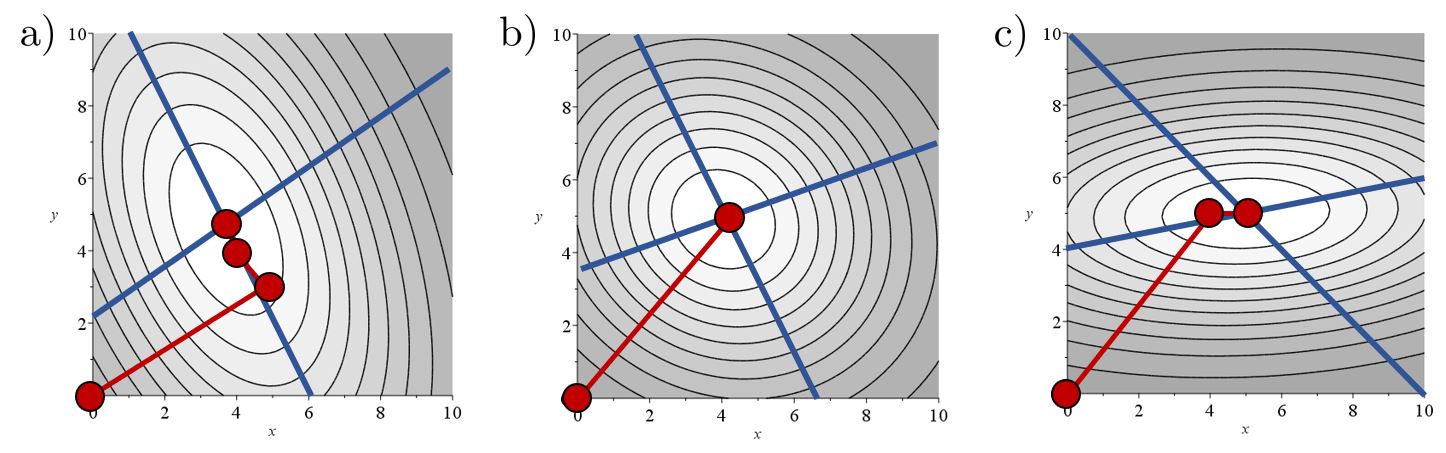}
\caption[tab]{Iterative convergence of the proposed model in 3 different systems.}
\label{fig:gaussour}
\end{figure}

Now, the convergence of the algorithm shown in equation \ref{ecuxbeta} is evaluated through simple examples. We present, in Figure \ref{fig:gaussour}, 3 cases of simple systems of equations where the Gauss-Seidel method is divergent and the proposed formulation, without updating and with $\beta=0$, converges within a few iterations. The cases are:\\

Case $a$:

\begin{equation}
-0.7\, x_1+x_2=2
\end{equation}
\begin{equation}
2\,x_1+x_2=12
\end{equation}

With eigenvalues $(-1.5,\;1.8)$ and eigenvectors $[(-0.78,\;-0.37),\;(0.62,\;-0.93)]$.\\

Case $b$:

\begin{equation}
-0.7\, x_1+2\,x_2=7
\end{equation}
\begin{equation}
2\,x_1+x_2=13
\end{equation}

With eigenvalues $(-2.0,\;2.3)$ and eigenvectors $[(-0.83,\;-0.55),\;(0.55,\;-0.83)]$.\\

Case $c$:

\begin{equation}
x_1-5\,x_2=-20
\end{equation}
\begin{equation}
2\,x_1+2\,x_2=20
\end{equation}

With eigenvalues $(1.5+3.1\,i,\; 1.5-3.1\,i)$ and eigenvectors $[(0.84,\;0.84),\;(-0.1-0.5\,i,\;-0.1+0.5\,i)]$.\\

And graphically:\\

In the system analysed, all the lines cut in the same point, which is the solution of the system in the $\mathds{R}^n$ space, where $n$ is the dimension of the system. This is shown in Figure \ref{fig:gaussour}, where if the system is square, the crossing point is always unique and its dimension is the range of the matrix $\mathcal{A}$ of the lineal system $\mathcal{A}\cdot \textbf{x}=\textbf{b}$. And this will be unique for any error introduced in the vector $\textbf{b}$. \\

\begin{figure}[H]
\centering
\includegraphics[width=0.4\textwidth]{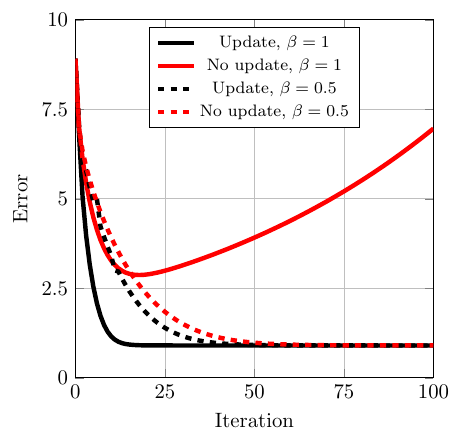}
\caption[tab]{Convergence of the proposed formulation under different approaches.}
\label{fig:noaquarelines}
\end{figure}

But in case of a non-square system like the one shown in Figure \ref{fig:noaquarelines} and equations \ref{matexamp} and \ref{bexamp}, the insertion of an error in the $\textbf{b}$ means that the crossing point is not unique. Instead, we will have a set of points which are crossings defined by different combinations of $n$ lines in the $\mathds{R}^n$ space, which defines the region where the minimum is located.  \\

Here the convergence is more difficult due to the non existence of a solution for the insertion of noise in the $\mathcal{A}$ matrix and its non-square size, so it is used for the study of the convergence of the different formulations presented in equations \ref{ecu1}, \ref{ecu2}, \ref{ecu3}, \ref{ecu4}, \ref{ecu5}. The error is calculated as the modulus of the $b$ vector obtained, and the one imposed in the equation $\mathcal{A}\cdot \textbf{x}= \textbf{b}$, where $u$ is the vector of unknowns that we are solving with the proposed algorithm.\\

The next case is a simple example without exact solution chosen to produce the divergence in the algorithm proposed, in the case without any updating, after a random evaluation of several matrices. In general, the divergence is not produced in most of the cases but this case is useful to present the error propagation in all the cases. The $\mathcal{A}$ matrix and the $\textbf{b}$ vector are:\\

\begin{equation}
\mathcal{A}=
\begin{bmatrix}
-8.11 & 2.75 & 9.52 & 6.57 & 1.17\\
6.35 & 9.21 & -7.61 & 8.51 & 9.91 \\
-7.43 & 1.12 & -0.64 & -8.75 & 4.12 \\
3.99 & 5.68 & -8.49 & 9.07 & -5.43 \\
6.00 & 5.33 & -9.56 & 1.74 & -5.62 \\
2.22 & -2.10 & -1.87 & -2.67 & 6.00 \\
-1.11 & 3.97 & 7.73 & 5.24 & 8.64 \\
7.70 & -4.45 & -2.38 & -9.23 & -2.75 \\
4.27 & -4.06 & -0.09 & -2.13 & -8.05 \\
0.72 & -0.53 & 8.69 & 1.02 & -6.85
\end{bmatrix}
\label{matexamp}
\end{equation}

\begin{equation}
\textbf{b}^T=
\begin{bmatrix}
-0.29 & -2.09 & 2.33 & 0.16 & 4.32 & -3.82 & -0.55 & 3.33 & 2.09 & 4.51
\end{bmatrix}
\label{bexamp}
\end{equation}

This is implemented in Figure \ref{fig:noaquarelines}:\\

This result shows the divergence of the implementation without the update and with $\beta=1$, i.e. where any error can be propagated freely. Any other case, where the errors are distributed and mitigated, make the method convergent to the solution, which is not $0$ because the matrix is non-square and there is some error inserted into the system of equations.

\subsection{Tensor formulation}

In this section we study the double dot product as an extension of the dot product studied in the system of equations. The high order tensors can be reduced through the application of the Voight's  notation \cite{Mnik2021} to lower order matrices, which implies as well a reduction in the order of the tensorial operation when the operation is computed using Voight's notation. It means that a general formulation can be achieved for the dot product between first order and second order tensors. By using the Voight's notation and the dot product the result can be extrapolated to high order tensors under double dot product. Within this premise, a general formulation is deducted from the same reasoning as the previous results.\\

The double dot product is an operation that reduces the dimension of the resulting tensor, when compared to the tensors involved in the operation. It means that a lower number of equations is available for the optimization of the tensor which is the solution of the equation and, as a general restriction, it can only be applied when the order of the resulting tensor $\mathbb{C}$ and the unknown tensor $\mathbb{B}$ are the same, perhaps the dimensions are not necessary equal as in the case of nested matrices, this makes the system compatible and feasible to be solved.\\

In the case of a high order tensor, we have as a general form $\mathbb{A}:\mathbb{X}=\mathbb{B}$, which in index form for 4 order tensor is $\mathbb{A}:\mathbb{X}=\mathbb{B}$, and the exponential kernel can be studied as: \\

\begin{equation}
\mathpzc{g}(\mathbb{X})=\exp \left[- \frac{(\mathbb{B}_{ijmn}-\mathbb{A}_{ijkl}\mathbb{X}_{klmn})^T(\mathbb{B}_{ijmn}-\mathbb{A}_{ijkl}\mathbb{X}_{klmn})}{\gamma} \right]  
\end{equation}

Which is the equivalent of equation \ref{ecu:g1} for the double dot product, and where:
%

\begin{equation}
\mathpzc{f}_{\hat{\textbf{k}}\hat{\textbf{l}}\hat{\textbf{m}}\hat{\textbf{n}}}(\mathbb{X})=\frac{\partial \mathpzc{g}(\mathbb{X})}{\partial x_{\hat{\textbf{k}}\hat{\textbf{l}}\hat{\textbf{m}}\hat{\textbf{n}}}}=0
\end{equation}

And, considering that for the double dot product in 4th order tensors only two indices are shared between $\mathbb{A}$ and $\mathbb{X}$, the component $\widehat{x}_{\hat{\textbf{k}}\hat{\textbf{l}}\hat{\textbf{m}}\hat{\textbf{n}}}$ can be expressed as:


\begin{equation}
\widehat{x}_{\hat{\textbf{k}}\hat{\textbf{l}}\hat{\textbf{m}}\hat{\textbf{n}}}=\frac{ \textbf{A}_{ij\hat{\textbf{k}}\hat{\textbf{l}}}^T \left(\textbf{B}_{ij\hat{\textbf{m}}\hat{\textbf{n}}}- \mathbb{A}_{ijop}\breve{\textbf{X}}_{op\hat{\textbf{m}}\hat{\textbf{n}}}(1-\delta_{o\hat{\textbf{k}}}\delta_{p\hat{\textbf{l}}})\right)}{\textbf{A}_{ij\hat{\textbf{k}}\hat{\textbf{l}}}^T \textbf{A}_{ij\hat{\textbf{k}}\hat{\textbf{l}}}}
\label{4orderdoubledot}
\end{equation}

\subsubsection{Generic iterative formulation}

The equation used previously for the calculation of the unknown tensor $X$ can be extended to the most common  operations with higher order tensors (single and double contraction) as $A\circ X=B$, where $\circ$ can be the dot or double dot operations. This generic equation \ref{generic} is the one used for all the results presented in this paper:

\begin{equation}
\widehat{x}_{\gamma}=(1-\beta)\tilde{x}_{\gamma}+\beta\frac{ A_{r_{\gamma}}^T \circ \left( B_{r_{\gamma}}-A\circ \breve{X}_{r_{\gamma}}+A_{r_{\gamma}}\tilde{x}_{\gamma}\right)}{A_{r_{\gamma}}^T\circ A_{r_{\gamma}}}
\label{generic}
\end{equation}

Where the operation $\circ$ is substituted by the dot or double dot product, and the reduced tensors are defined as:

\begin{itemize}
\item[$\rightarrow$] $\beta$: weight parameter of the solution.
\item[$\rightarrow$] $\gamma$: set of indices that defines the position of the variable $x_{\gamma}$ studied, which is a component of the high order unknown tensor $X$.
\item[$\rightarrow$] $\widehat{x}_{\gamma}$: value of the variable studied at time $t+1$.
\item[$\rightarrow$] $\tilde{x}_{\gamma}$: value of the variable studied at time $t$.
\item[$\rightarrow$] $A_{r_{\gamma}}$: tensor $A$, of any order, reduced by blocking the indices defined by $\gamma$.
\item[$\rightarrow$] $B_{r_{\gamma}}$: tensor $B$, of any order, reduced by blocking the indices defined by $\gamma$.
\item[$\rightarrow$] $\breve{X}_{r_{\gamma}}$: tensor $X$, of any order, reduced by blocking the indices defined by $\gamma$.
\end{itemize}

To ensure that the system can be solved, it needs to be guaranteed that the operation $A_{r_{\gamma}}^T\circ A_{r_{\gamma}}$ has a scalar as a result, which determines the size of $A$ that can be used for every operation.\\

For example, the equation \ref{4orderdoubledot} can be expressed as:\\

\begin{equation}
\widehat{x}_{\hat{\textbf{k}}\hat{\textbf{l}}\hat{\textbf{m}}\hat{\textbf{n}}}=\frac{ \mathbb{A}_{r_{\hat{\textbf{k}}\hat{\textbf{l}}}}^T: \left(\mathbb{B}_{r_{\hat{\textbf{m}}\hat{\textbf{n}}}}- \mathbb{A}:\breve{\mathbb{X}}_{r_{\hat{\textbf{m}}\hat{\textbf{n}}}}(1-\delta_{o\hat{\textbf{k}}}\delta_{p\hat{\textbf{l}}})\right)}{\mathbb{A}_{r_{\hat{\textbf{k}}\hat{\textbf{l}}}}^T : \mathbb{A}_{r_{\hat{\textbf{k}}\hat{\textbf{l}}}}}
\end{equation}

Which is an implementation of the proposed notation for the reduction of tensors by blocking some indices.

\subsection{Non square systems without exact solution}

\begin{figure}[H]
\centering
\includegraphics[width=0.8\textwidth]{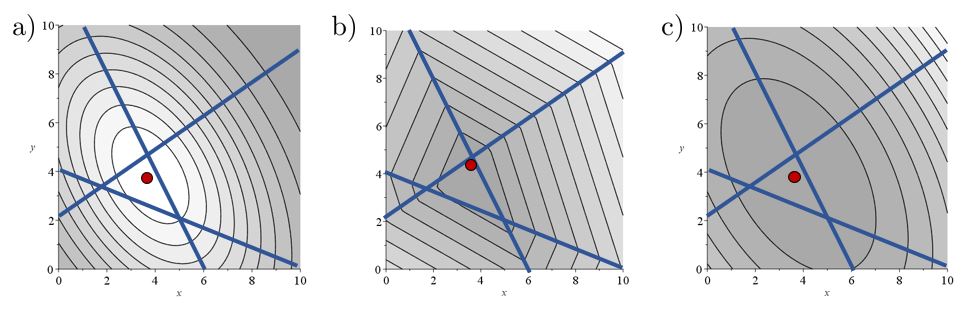}
\caption[tab]{Maximum in red and the three lines that represent the equations of the system. Exponential kernel with $\gamma=2$ (a), linear function (b) and quadratic (c).}
\label{fig:centre}
\end{figure}

A typical case of application is a noisy dataset of linear equations in non-square systems. In those, there are more equations than variables, but the induced noise means that there is no exact solution for the system, and the problem turns into a minimization process to find the point that reduces the objective function imposed. And here, there is a difference between the solutions depending on the objective function chosen and the solver used. For example, Figure \ref{fig:centre} shows a system composed by the equations:\\

\begin{equation}
-0.7\, x_1+x_2=2
\end{equation}
\begin{equation}
2\,x_1+x_2=12
\end{equation}
\begin{equation}
0.4\, x_1 +x_2=4
\end{equation}

Which are represented in the Figure \ref{fig:centre} by the blue lines, and the contour plot in gray scale shows the value of the objective function used. Figure \ref{fig:centre}a shows the exponential kernel presented in equation \ref{ker_obj}, Figure \ref{fig:centre}b shows lineal objective function defined in equation \ref{lineal_obj}, and Figure 3c the quadratic objective function defined by equation \ref{cuad_obj}. 

\begin{equation}
L_{obj}(x_1,x_2)=|-0.7\, x_1+x_2-2|+|2\,x_1+x_2-12|+|0.4\, x_1 +x_2-4|
\label{lineal_obj}
\end{equation}

\begin{equation}
Q_{obj}(x_1,x_2)=(-0.7\, x_1+x_2-2)^2+(2\,x_1+x_2-12)^2+(0.4\, x_1 +x_2-4)^2
\label{cuad_obj}
\end{equation}

\begin{equation}
K_{obj}(x_1,x_2)=\exp\left[ -\left( \frac{(-0.7\, x_1+x_2-2)^2+(2\,x_1+x_2-12)^2+(0.4\, x_1 +x_2-4)^2}{\gamma}\right) \right] 
\label{ker_obj}
\end{equation}

And so, the quadratic and exponential kernel objective functions are related as:\\

\begin{equation}
K_{obj}(x_1,x_2)=\exp\left[-\frac{Q_{obj}(x_1,x_2)}{\gamma} \right]
\label{ker2_obj}
\end{equation}

Here, the minimum of the quadratic objective function is the maximum of the exponential kernel. In those cases without a minimum equal to 0 due to the error in the solution and the non-square system, and its consequent lack of a exact solution, the function $Q_{obj}(x_1,x_2)$ does not have a kernel, and thus the $K_{obj}(x_1,x_2)$ function can not satisfy $K_{obj}(x_1,x_2)=1$ in any point of the space. But the exponential kernel ensures that the minimum of $Q_{obj}(x_1,x_2)$ is exactly the maximum of $K_{obj}(x_1,x_2)$. For this reason, even if there is no kernel, we will use the exponential kernel theory to solve the non-square noisy systems of equations.\\

Since there is no exact solution, the lines in Figure \ref{fig:centre} are not crossing in a single point. And instead, a set of crossings between subsets of lines is generated, which defines the boundary of a region that holds the optimum of the full problem. For example, if we have a system with $m$ equations in a space of dimension $n$ (with $m>n$), without exact solution, and there will be subsets of $n$ equations that will have an exact solution and will define a crossing point in the space of dimension $n$.\\

More formally, in a system with an exact solution all the lines cut in the same point, which is the solution of the system in the $\mathds{R}^n$ space, where $n$ is the dimension of the system. This is shown in Figure \ref{fig:gaussour}, where if the system is square, always the crossing point is unique and its dimension is the range of the matrix $\mathcal{A}$ of the lineal system $\mathcal{A}\cdot \textbf{u}=\textbf{b}$. And this will be unique for any error introduced in the vector $\textbf{b}$. But in the case of a non-square system like the one shown in Figure \ref{fig:noaquarelines}, the insertion of an error in the $\textbf{b}$ means that the crossing point is not unique. Instead, we will have a set of points which are crossings defined by different combinations of $n$ lines in the $\mathds{R}^n$ space. In a system with a $\mathcal{A}$ matrix of $n\times m$, each center will be defined for every submatrix of $n\times n$ within the system. We define $\kappa$ as this number of centres is determined by the combinatorics theory as:\\

\begin{equation}
\kappa={m \choose n}= \frac{m!}{n!\left( m-n\right)! }
\end{equation}

Which, for this particular case with 3 equations and 2 variables, is:

\begin{equation}
\kappa={3 \choose 2}= \frac{3!}{2!\left( 3-2\right)! }=\frac{6}{2}=3
\end{equation}

Those are the 3 crossings between the lines shown in Figure \ref{fig:centre}, where the maximum or minimum of each objective function is represented with a red dot. It can be seen that the the lineal optimum differs from the one obtained by both the quadratic and the exponential kernel, which are coincident.

In Figure \ref{fig:centre}, it can be clearly seen that the 3 crossing points of the lines, which are the centres, define a convex region where the optimum is located, with independence of the objective function. Also, the centre of the linear problem differs from the one calculated by both the quadratic and the exponential kernel.

\subsection{Algorithmic complexity for the inversion of matrices}

In this section we will analyse the number of operations needed to calculate an inverse or pseudoinverse of a matrix. For this, we start with the equation for the operation $\mathcal{A} \cdot \mathcal{X} =\mathcal{B}$ where $\mathcal{A}$ is the known matrix that we want to invert, $\mathcal{X}$ is the inverse or pseudoinverse, and $\mathcal{B}$ is the identity matrix $\mathcal{I}$. So, in index notation:

\begin{equation}
\widehat{\mathcal{X}}_{\hat{\textbf{j}}\hat{\textbf{l}}}=(1-\beta)\tilde{\mathcal{X}}_{\hat{\textbf{j}}\hat{\textbf{l}}}+\beta\frac{  \mathcal{A}_{k\hat{\textbf{j}}}^T\left( \mathcal{I}_{k\hat{\textbf{l}}}-\mathcal{A}_{ks} \breve{\mathcal{X}}_{s\hat{\textbf{l}}}+\mathcal{A}_{k\hat{\textbf{j}}} \tilde{\mathcal{X}}_{\hat{\textbf{j}}\hat{\textbf{l}}}\right)}{ \mathcal{A}_{k\hat{\textbf{j}}}^T \mathcal{A}_{k\hat{\textbf{j}}}}
\label{ecupseudo}
\end{equation}

The capabilities of this equation, compared with the Moore-Penrose pseudoinverse is studied later on in this paper.

In this section we are going to calculate the number of operations needed to carry out the calculation. For this, the equation \ref{ecupseudo} is used to calculate one iteration of the algorithm. We assume that the matrix $\mathcal{A}$ has a size of $m \times n$, $\mathcal{X}$ a size of $n \times m$ and $\mathcal{I}$ a size of $m \times m$. Here we can crop this equation into the different computations, considering as well the number of operations needed for a full iteration, considering the range of values for the blocked indices of the operation:

\begin{itemize}
\item[$\rightarrow$] $\mathcal{A}_{k\hat{\textbf{j}}}^T \mathcal{I}_{k\hat{\textbf{l}}}$: this operation is the dot product of two vectors of size $m$. This has $m$ multiplications and $m-1$ sums, so the total number of operations is $2m-1$. This is calculated $n$ times in one full iteration, so the total number of operations is $2mn-n$.

\item[$\rightarrow$] $-\mathcal{A}_{k\hat{\textbf{j}}}^T \left( \mathcal{A}_{ks} \breve{\mathcal{X}}_{s\hat{\textbf{l}}}\right)$: inside the parenthesis we have a dot product between a matrix and a vector, resulting in a vector of size $m$. This operation is the evaluation of $m$ dot products between vectors of dimension $n$, so we have $m(2n-1)$ operations, and to do a full iteration we need to calculate this $m$ times, so we have $2nm^2-m^2$ operations. Now, to operate with the term outside the parenthesis, we have a dot product between two vectors, which again is $2m-1$ operations, carried out $n$ times for a full operation with a total of $2mn-n$ operations. So, the equation of this point is solved for 1 complete iteration with $2mn-n+2nm^2-m^2$ operations.

\item[$\rightarrow$] $\left( \mathcal{A}_{k\hat{\textbf{j}}}^T \mathcal{A}_{k\hat{\textbf{j}}}\right)  \tilde{\mathcal{X}}_{\hat{\textbf{j}}\hat{\textbf{l}}}$: this is a dot product between two vectors of dimension $m$ (inside the parenthesis), multiplied by an scalar. So, inside the parenthesis there are $2m-1$ operations, which are repeated $n$ times to make the full iteration (i.e. $2mn-n$ operations). Then, the product by an scalar is another operation that is done $mn$ times, considering the dimensions of the fixed indices that are repeated to calculate the different coefficients. In total, this part of the equation involves $2mn-n+mn+1$ operations.

\item[$\rightarrow$] $\mathcal{A}_{k\hat{\textbf{j}}}^T \mathcal{A}_{k\hat{\textbf{j}}}$: this was studied before with a total number of $2mn-n$ operations.  
\end{itemize}

To this, the final division of equation \ref{ecupseudo} adds another operation that is done $mn$ times for each iteration, as well as the multiplication by $\beta$ and the sum. So, for every iteration in all the indices we have a total of $s(2nm^2-m^2+9mn-3m)$ operations for all the iterations $s$ needed for the convergence. So, by using the $\mathcal{O}$ notation, the complexity of the algorithm to calculate the pseudoinverse of a matrix is of the order of $\mathcal{O}(nm^2)$. This is exactly the algorithmic complexity of the Singular Value Decomposition (SVD) \cite{SVD} used to calculate the pseudoinverse with the Moore-Penrose method \cite{penrose_1955}.\\

 In an square system, where $m=n$, the complexity of the proposed algorithm is $\mathcal{O}(n^3)$, and this algorithm can be compared with other methodologies used for the matrix inversion. The table below show the complexity and capabilities of some of the most used algorithms:

\begin{table}[H]
\centering
\begin{tabular}{ c c c c c c}
 Algorithm & G-J & St & oC-W & M-P & Proposed \\ 
Complexity & $\mathcal{O}(n^3)$ & $\mathcal{O}(n^{2.807})$ & $\mathcal{O}(n^{2.373})$ & $\mathcal{O}(n^3)$ &  $\mathcal{O}(n^3)$ \\  
 Iterative & Yes & Yes & Yes & Yes & Yes \\
 Preconditioning & Yes & Yes & Yes & No & No\\
 Jacobian required & No & No & No & Yes & No\\
 Suitable for pseudoinverse & No & No & No & Yes & Yes  
\end{tabular}
 \caption{Comparison of some of the most used algorithms for matrix inversion and the characteristics of the proposed model in this paper.}
 \label{table_alg}
 \end{table}

The Table \ref{table_alg} shows a comparison between some algorithms and the one proposed in this work. Those models are denoted by their initials, where G-J is the Gauss-Jordan reduction, St is the Strassen algorithm \cite{Strassen1969}, oC-W is the optimized Coppersmith-Winograd algorithm \cite{Williams2012}, and M-P is the Moore-Penrose pseudoinverse. The second row of the table shows the big-$\mathcal{O}$ complexity of each one, and then if any preconditioning is needed. It means that some arrangements or decision needs to be considered in order to make the calculation, such as the reordering of the matrix, the selection of sub-matrices or the order of the operations. Then, only the Moore-Penrose method needs to calculate the Jacobian to reach the solution, but except ours, any other of the analysed can calculate a pseudoinverse in a non-square system without the use of the Jacobian.\\

For the complexity, all of them are close and the proposed algorithm has the same complexity as the Moore-Penrose, which is the closest in terms of capabilities. 

\section{Results}

In this section we present different results of the algorithm. First, we present an analytical reduction for small systems, which lies on the Cramer formulation \cite{Kosinski2001}. Then, this is generalized for a system of equations to be solved iteratively, and after a simplification for diagonal dominant systems we reach the Jacobi and Gauss-Seidel formulations. Then the general formulation is compared with the Moore-Penrose pseudoinverse for a different range of high order tensors and nested non-square matrices without exact solution. Reaching exactly the same results linearly with the algorithm proposed, instead of by using the least square method implemented in the Moore-Penrose pseudoinverse. \\

\subsection{Cramer}

Consider a system composed by 2 equations and 2 unknowns ($x_1$, $x_2$), and the equation \ref{ecuxsyslin} that extract one variable as a function of the others. If we apply this equation to $x_1$ and $x_2$ we have:\\

\begin{equation}
x_1=\frac{b_1a_{11}-a_{11}a_{12}x_2+b_2a_{21}-a_{21}a_{22}x_2 }{a_{11}^2+a_{21}^2}
\label{equ_x1}
\end{equation}

\begin{equation}
x_2=\frac{b_1a_{12}-a_{12}a_{11}x_1+b_2a_{22}-a_{22}a_{21}x_1 }{a_{12}^2+a_{22}^2}
\label{equ_x2}
\end{equation}

If we use the results of equations \ref{equ_x1} and \ref{equ_x2} and solve the values of the variables as a function of the terms $a_{ij}$ and $b_i$, we obtain that:\\

\begin{equation}
x_1=\frac{a_{22}b_1-a_{12}b_2}{a_{11}a_{22}-a_{12}a_{21}}=\frac{
  \begin{vmatrix} b_1 & a_{12} \\ b_2 & a_{22} \end{vmatrix}  
}{\begin{vmatrix} a_{11} & a_{12} \\ a_{22} & a_{22} \end{vmatrix} }
\label{x1c}
\end{equation}

\begin{equation}
x_2=\frac{a_{11}b_2-a_{21}b_1}{a_{11}a_{22}-a_{12}a_{21}}=\frac{
  \begin{vmatrix}  a_{11} & b_1 \\ a_{21} & b_2 \end{vmatrix}  
}{\begin{vmatrix} a_{11} & a_{12} \\ a_{22} & a_{22} \end{vmatrix} }
\label{x2c}
\end{equation}

Where the equations \ref{x1c} and \ref{x2c} are the classical Cramer solutions for systems of equations. The formulation proposed has the same exact solution as the one presented by Cramer in 1750, and can be extended for a larger number of unknowns. This has the same potentialities and drawbacks as the Cramer´s rule, with a very high computational cost for large systems of equations. But the formulation presented in the equation \ref{ecuxsyslin} allows us to reach an iterative solution of the problem, which relies on the smoothness of the exponential function and will be studied in the next section.

\subsection{Gauss-Seidel}

The Gauss-Seidel formulation is an extension of the Jacobi algorithm, both widely used in the iterative calculation of square, symmetric, and diagonal dominant systems of equations.

The formulation deducted by Jacobi relies on a diagonal dominant, positive and symmetrical matrix $\mathcal{A}$, which are essential requisites for its convergence. The general formulation of the proposed algorithm without updating is:

\begin{equation}
x_j^{t+1}=\frac{\sum_i \left(b_i a_{ij}- a_{ij}\sum_{k \neq j}\left( a_{ik}x_k^{t}\right) \right) }{\sum_i a_{ij}^2}
\label{equ_xja}
\end{equation}

The equation \ref{equ_xja} proposed, can be simplified considering a diagonal dominant matrix, i.e. assuming that the terms of the diagonal are much bigger than the terms outside it. So, we are only going to consider, from equation \ref{equ_xja}, the terms of the sum where $i=j$, which is shown in equation \ref{equ_xj_jac}:

\begin{equation}
x_j=\frac{\left(b_j a_{jj}- a_{jj}\sum_{k \neq j}\left( a_{jk}x_k^{t}\right) \right) }{a_{jj}^2}\thickapprox\frac{\left(b_j -\sum_{k \neq j}\left( a_{jk}x_k^{t}\right) \right) }{a_{jj}}
\label{equ_xj_jac}
\end{equation}

Where the equation \ref{equ_xj_jac} is exactly the solution of Jacobi, which is a convergence of the proposed algorithm when the matrix $\mathcal{A}$ is assumed to be strongly diagonal dominant. It proves that the Jacobi formulation is a simplification of the proposed algorithm.

The Gauss-Seidel method is the same as Jacobi but with a gradual update of the variables. In this case, the algorithm proposed with the update of the variables is:
 
\begin{equation}
x_j^{t+1}=\frac{\sum_i \left(b_i a_{ij}- a_{ij}\sum_{k < j}\left( a_{ik}x_k^{t+1}\right)- a_{ij}\sum_{k > j}\left( a_{ik}x_k^{t}\right) \right) }{\sum_i a_{ij}^2}
\label{equ_xjaup}
\end{equation}

Imposing to the equation \ref{equ_xjaup} the same conditions used to deduct the Jacobi equation, i.e. assume that the matrix is strong diagonal dominant, we have:

\begin{equation}
x_j=\frac{\left(b_j a_{jj}- a_{jj}\sum_{k < j}\left( a_{jk}x_k^{t+1}\right)- a_{jj}\sum_{k > j}\left( a_{jk}x_k^{t}\right) \right) }{a_{jj}^2}\thickapprox \frac{\left(b_j- \sum_{k < j}\left( a_{jk}x_k^{t+1}\right)- \sum_{k > j}\left( a_{jk}x_k^{t}\right) \right) }{a_{jj}}
\label{equ_xj_gau_sei}
\end{equation}

Those results prove that both the Jacobi and Gauss-Seidel formulations are simplifications of the proposed algorithm for the specific particular case of a strong diagonal dominant square matrix. This simplification, from the global formulation proposed, is also an explanation for the divergence of the Jacobi and Gauss-Seidel methods, when the conditions of strong diagonal dominance are not fulfilled by the matrix. 

\subsection{Matrix and tensor calculation}

Here the application of the algorithm will be tested numerically in systems where the traditional linear solvers, such as the Gauss-Seidel, diverge. Those are high order tensors and non-square nested matrices with noise. In this section we prove that the proposed algorithm can reach the same capabilities as the Moore-Penrose pseudoinverse but with a linear approximation, instead of the quadratic approach of the least-square method. Also, there is no need for gradient calculation or square roots solvers, since the algorithm is linear and exactly the same for every step. The four cases presented will implement the algorithm in different operations and orders, both for tensors and non-square nested matrices. Those are compared with the Moore-Penrose results, even when there is no exact solution due to the artificially inserted noise in non-square systems, and the solution is the result of a minimization.

\subsubsection{Case 1: $\mathbf{A} \cdot \mathbf{X} =\mathbf{B} \; \rightarrow \;\mathcal{A}_{ij} \textbf{x}_j=\textbf{b}_i$}

This first case is the matrix formulation of a system of linear equations. The algorithm proposed, for this case will be:

\begin{equation}
\widehat{x}_\alpha=\frac{ \mathcal{A}_{\alpha}^T \cdot \left( \textbf{b}_{\alpha}-\mathcal{A}\cdot \breve{\textbf{x}}_{\alpha}+\mathcal{A}_{\alpha}\tilde{x}_\alpha\right)}{\mathcal{A}_{\alpha}^T \cdot \mathcal{A}_{\alpha}}
\end{equation}

Which in indicial notation is:

\begin{equation}
\widehat{x}_{\hat{\textbf{j}}}=(1-\beta)\tilde{x}_{\hat{\textbf{j}}}+\beta\frac{\mathcal{A}_{l\hat{\textbf{j}}}^T\left( \textbf{b}_{l}-\mathcal{A}_{lk}\breve{\textbf{x}}_{k}+\mathcal{A}_{l\hat{\textbf{j}}}\tilde{x}_{\hat{\textbf{j}}}\right)}{ \mathcal{A}_{l\hat{\textbf{j}}}^T \mathcal{A}_{l\hat{\textbf{j}}}}
\end{equation}

\begin{figure}[H]
\centering
\subfloat{\includegraphics[width=0.32\textwidth]{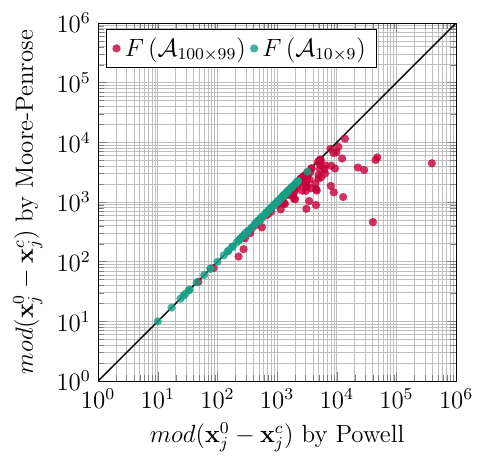}}
\hspace{0.1cm}
\subfloat{\includegraphics[width=0.32\textwidth]{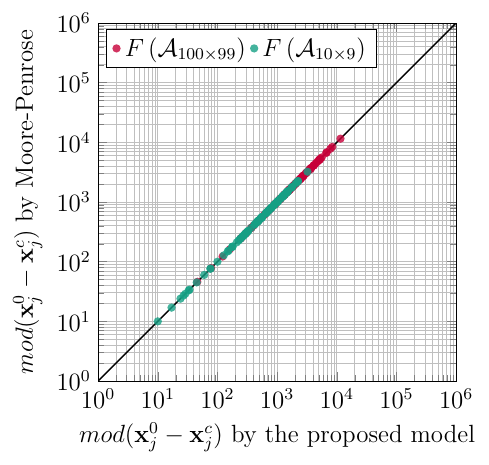}}
\hspace{0.1cm}
\subfloat{\includegraphics[width=0.32\textwidth]{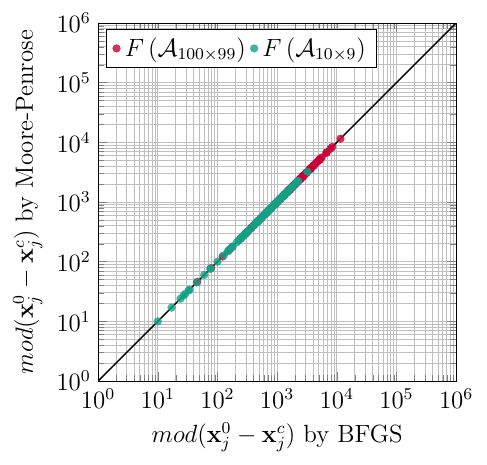}}
\caption{Comparison of different linear and quadratic algorithms with the Moore-Penrose pseudoinverse.}
\label{fig:2x1}
\end{figure}

This formulation will be used here to calculate the error considering random non-square matrices without exact solution. In the case of Figure \ref{fig:2x1} 100 random systems of 100 variables and 99 equations, and another 100 variables of 10 variables and 9 equations. In all the cases some noise were inserted in the $\mathbf{b}$ vector to ensure that there is no exact solution and the problem relies in a minimization process.\\

Figure \ref{fig:2x1} shows the error of different methodologies in those matrices. This error is measured as the modulus of the original $\mathbf{b}$ vector of the system minus the one resulting when the convergence is reached by every methodology. So, Figure \ref{fig:2x1} left shows the comparison between the results calculated with the Moore-Penrose pseudoinverse and the Powell method. The Powell method is linear and gives higher errors than the Moore-Penrose, which is based on the quadratic SLSQP algorithm.\\

Figure \ref{fig:2x1} centre shows the calculated result of the Moore-Penrose pseudoinverse compared with the proposed algorithm. It is clear that the correlation of the proposed model with the quadratic methods is exact.\\

\begin{figure}[H]
\centering
\subfloat{\includegraphics[width=0.32\textwidth]{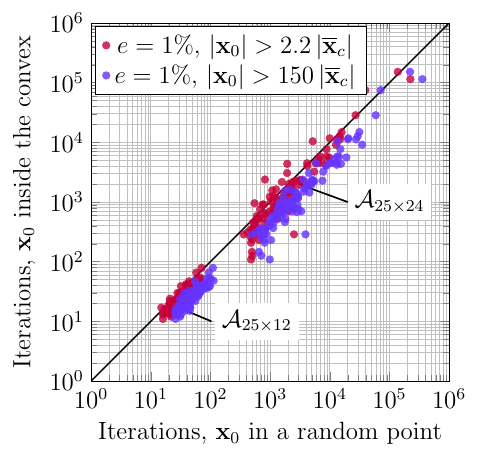}}
\hspace{0.1cm}
\subfloat{\includegraphics[width=0.32\textwidth]{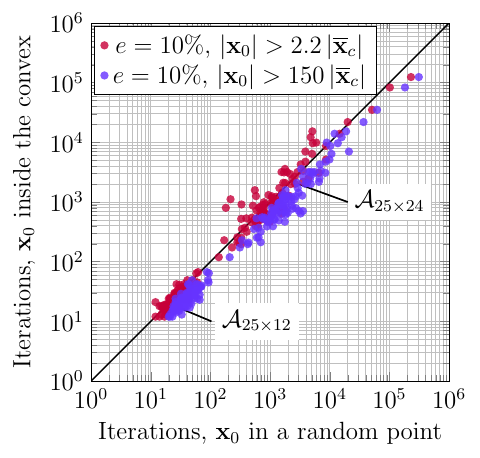}}
\caption{Iterations until convergence for different locations of the initial point.}
\label{fig:iterations}
\end{figure}

Since the method proposed is iterative, it needs to be studied the influence of the starting point of the iterations. As previously presented, there is a set of centres in every non-square matrices without exact solution. Each one is a partial solution of the system and defines a convex where the minimization of the system of equations is located. So, Figure \ref{fig:iterations} shows a comparison between the iterations needed until convergence when the starting point is inside the convex, compared with the iterations needed with the initial point located outside the convex. Two cases of initialization are studied, one closer to the convex than the other (i.e. 2.2 and 150 times the centroid of the centres). 800 random cases with two non-square matrix sizes and an error inserted to avoid the exact solution.\\

All the cases presented in Figure \ref{fig:iterations} reach exactly the same solution, but with a different number of iterations. And it probes that the starting point of the iterations only affects to the computational cost to reach the solution, but the convergence to the same point is always reached for every starting point.

\subsubsection{Case 2: $\mathbf{A} \cdot \mathbf{X} =\mathbf{B} \; \rightarrow \; \mathcal{A}_{ij}\mathcal{X}_{jl}=\mathcal{B}_{il}$}

Here we study a system of equation involving a set of unknowns $\mathbf{X}$ of order 2, operated with $\mathbf{A}$ and $\mathbf{B}$, of order two as well. The equation is:

\begin{equation}
\widehat{x}_\alpha=\frac{ \mathcal{A}_{\alpha}^T \cdot \left( \mathcal{B}_{\alpha}-\mathcal{A}\cdot \breve{\mathcal{X}}_{\alpha}+\mathcal{A}_{\alpha}\tilde{x}_\alpha\right)}{\mathcal{A}_{\alpha}^T\cdot \mathcal{A}_{\alpha}}
\end{equation}

And in indicial form:

\begin{equation}
\widehat{\mathcal{X}}_{\hat{\textbf{j}}\hat{\textbf{l}}}=(1-\beta)\tilde{\mathcal{X}}_{\hat{\textbf{j}}\hat{\textbf{l}}}+\beta\frac{  \mathcal{A}_{k\hat{\textbf{j}}}^T\left( \mathcal{B}_{k\hat{\textbf{l}}}-\mathcal{A}_{ks} \breve{\mathcal{X}}_{s\hat{\textbf{l}}}+\mathcal{A}_{k\hat{\textbf{j}}} \tilde{\mathcal{X}}_{\hat{\textbf{j}}\hat{\textbf{l}}}\right)}{ \mathcal{A}_{k\hat{\textbf{j}}}^T \mathcal{A}_{k\hat{\textbf{j}}}}
\label{equ:2x2=2}
\end{equation}

Equation \ref{equ:2x2=2} is applied to different cases, to calculate directly the Moore-Penrose pseudoinverse and to solve the unknowns $\mathbf{X}$.\\

\textit{Calculation of $\mathcal{X}$ in a non-square system without exact solution}\\

\begin{figure}[H]
\centering
\subfloat{\includegraphics[width=0.32\textwidth]{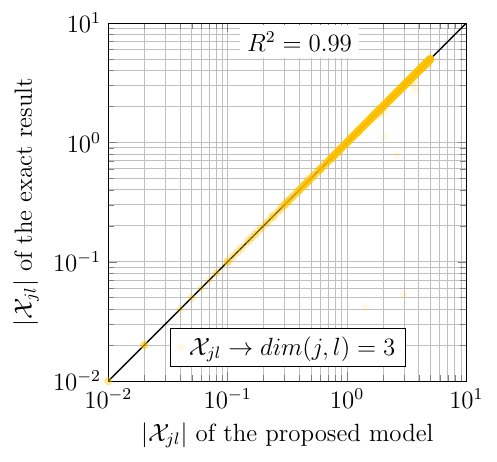}}
\subfloat{\includegraphics[width=0.32\textwidth]{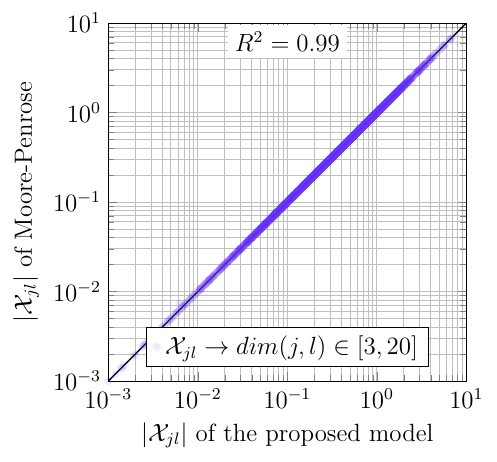}}
\subfloat{\includegraphics[width=0.32\textwidth]{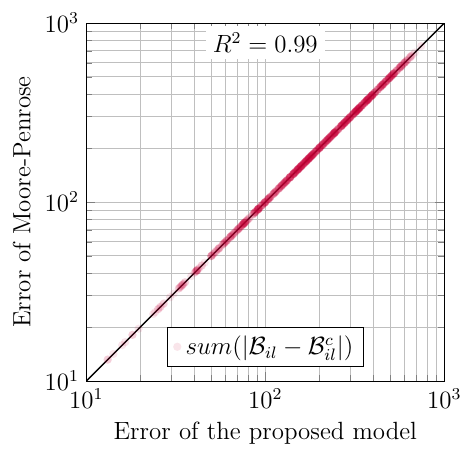}}
\caption{Comparison of the proposed model with the Moore-Penrose pseudoinverse for a tensor (left), and a nested matrix (centre) with its error (right).}
\label{fig:comparison2x2=2}
\end{figure}

Figure \ref{fig:comparison2x2=2} left shows the application of equation \ref{equ:2x2=2} to calculate tridimensional second order tensors, given random tensors $\mathbf{A}$ and $\mathbf{B}$ of the same order and dimension. The plot compares component by component 200 different calculated tensors with its exact counterpart, which are 1800 data points in total. Figure \ref{fig:comparison2x2=2} centre shows the same result but with non-square random matrices without exact solution. Over 200 cases, and since the plot compares component by component the resulting matrix, it has around 100000 data points in the validation against the components of the Moore-Penrose pseudoinverse. Figure \ref{fig:comparison2x2=2} right shows the error of the proposed model versus the error of the Moore-Penrose pseudoinverse, where this error is calculated as the modulus of the difference between the random $\mathcal{B}$ imposed in the system and the one calculated as $\mathcal{A} \cdot \mathcal{X}$, where $\mathcal{X}$ is the solution of the system calculated.

It shows the total agreement between the algorithm proposed and the Moore-Penrose solution for a wide variety of random systems.\\

\textit{Calculation of the Moore-Penrose pseudoinverse of $\mathcal{A}$}\\

The previous results lead us to use the algorithm proposed for the calculation of the direct calculation of the Moore-Penrose pseudoinverse, instead of comparing the global solution as in Figure \ref{fig:comparison2x2=2}. This is simply done by imposing that the system is equal to the identity matrix as: $A_{ij}X_{ji}\approx I_{ii}$, so the unknowns is exactly the Moore-Penrose pseudoinverse. In this case, we use the equation \ref{ecupseudo}.

\begin{figure}[H]
\centering
\subfloat{\includegraphics[width=0.32\textwidth]{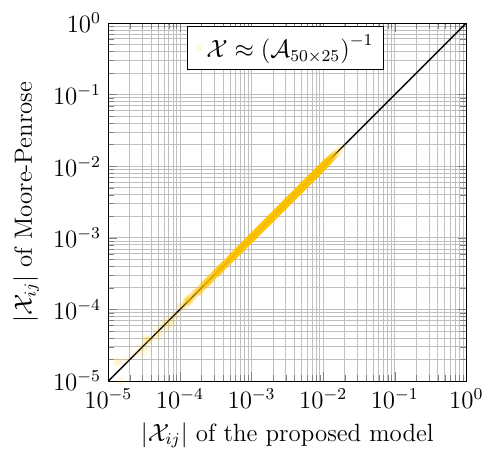}}
\hspace{0.1cm}
\subfloat{\includegraphics[width=0.32\textwidth]{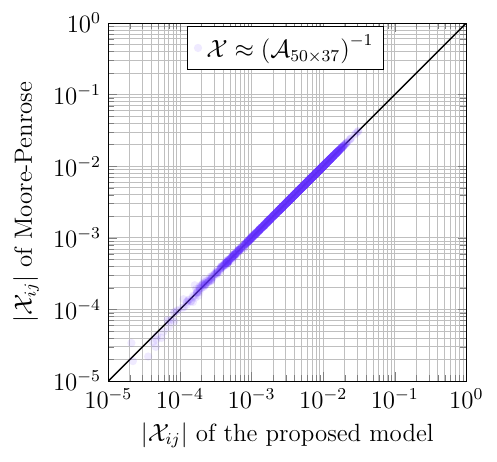}}
\hspace{0.1cm}
\subfloat{\includegraphics[width=0.32\textwidth]{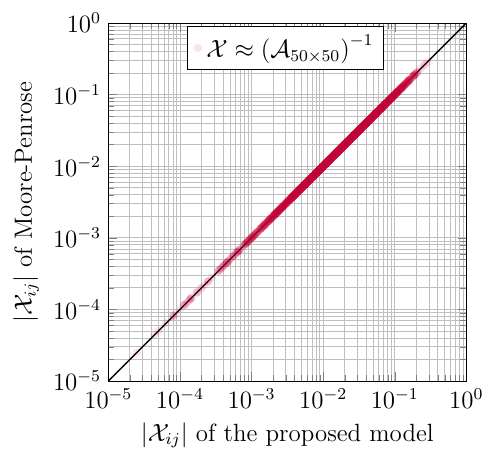}}\\

\subfloat{\includegraphics[width=0.32\textwidth]{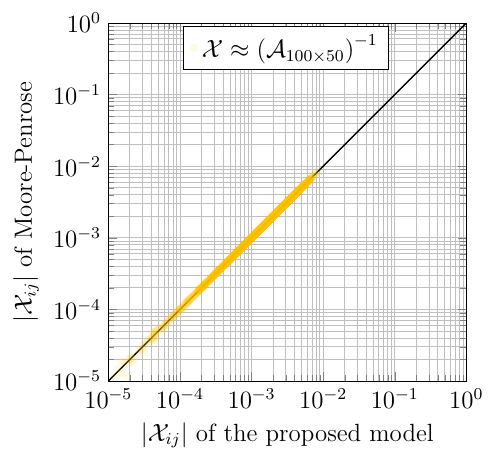}}
\hspace{0.1cm}
\subfloat{\includegraphics[width=0.32\textwidth]{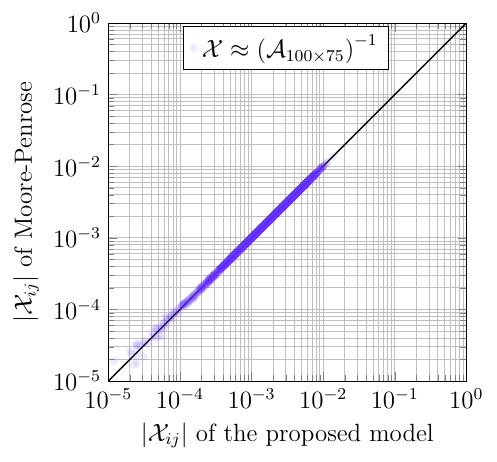}}
\hspace{0.1cm}
\subfloat{\includegraphics[width=0.32\textwidth]{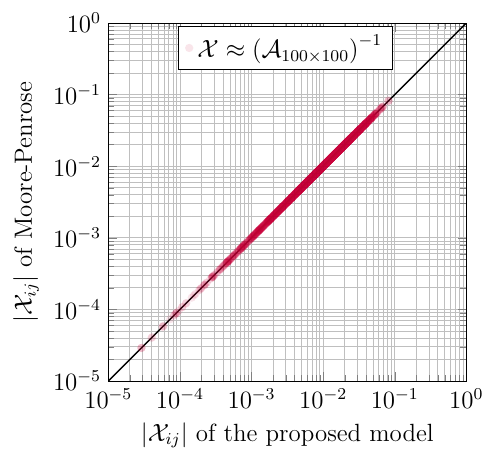}}
\caption{Comparison of the coefficients of the Moore-Penrose pseudoinverse with the proposed algorithm. Validation with approximately 1.5 million of coefficients of the pseudoinverses of 300 random matrices with different sizes.}
\label{fig:MPcomparison}
\end{figure}

Figure \ref{fig:MPcomparison} shows the comparison, coefficient by coefficient, between the result obtained by the proposed algorithm and the Moore-Penrose pseudoinverse. All are calculated in 300 completely random non-square matrices of different sizes, which is in total more than 1.5 millions of data points with a $R^2$ above $0.999$. 

The results of Figure \ref{fig:MPcomparison} are very interesting since the Moore-Penrose pseudoinverse is based on an iterative quadratic logarithm, and the proposed equation is an iterative linear logarithm. 

\subsubsection{Case 3: $\mathbf{A} : \mathbf{X} =\mathbf{B} \; \rightarrow \; \mathbb{A}_{ijkl}\mathcal{X}_{kl}=\mathcal{B}_{ij}$}

In this case, the order of $\mathbf{A}$ is equal to four, and the operation studied is the double dot product.

\begin{equation}
\widehat{x}_\alpha=\frac{ \mathbb{A}_{\alpha}^T : \left( \mathcal{B}_{\alpha}-\mathbb{A}: \breve{\mathcal{X}}_{\alpha}+\mathbb{A}_{\alpha}\tilde{x}_\alpha\right)}{\mathbb{A}_{\alpha}^T: \mathbb{A}_{\alpha}}
\end{equation}

$\breve{\mathcal{X}}_{\alpha}^{t}$ is composed by all the components $\alpha$ at time $t+1$. Where $\alpha$ is the element of the tensor $\mathcal{B}$ evaluated. In index form is:\\

\begin{equation}
\widehat{\mathcal{X}}_{\hat{\textbf{k}}\hat{\textbf{l}}}=(1-\beta)\tilde{\mathcal{X}}_{\hat{\textbf{k}}\hat{\textbf{l}}}+\beta\frac{ \mathbb{A}_{ij\hat{\textbf{k}}\hat{\textbf{l}}}^T\left( \mathcal{B}_{ij}-\mathbb{A}_{ijop}\breve{\mathcal{X}}_{op}+\mathbb{A}_{ij\hat{\textbf{k}}\hat{\textbf{l}}}\tilde{\mathcal{X}}_{\hat{\textbf{k}}\hat{\textbf{l}}}\right)}{\mathbb{A}_{ij\hat{\textbf{k}}\hat{\textbf{l}}}^T \mathbb{A}_{ij\hat{\textbf{k}}\hat{\textbf{l}}}}
\end{equation}

\begin{figure}[H]
\centering
\subfloat{\includegraphics[width=0.31\textwidth]{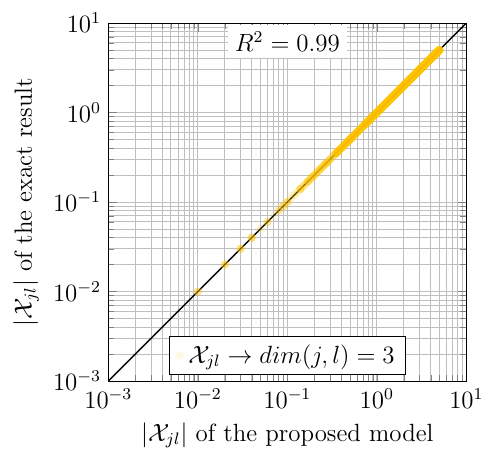}}
\subfloat{\includegraphics[width=0.31\textwidth]{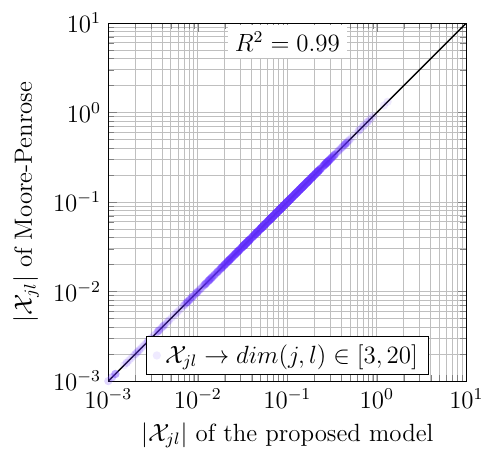}}
\subfloat{\includegraphics[width=0.30\textwidth]{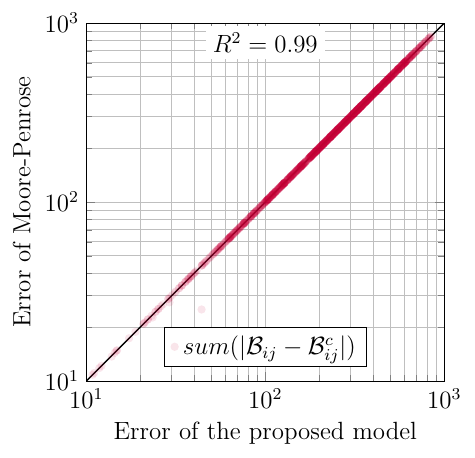}}
\caption{Comparison of the Moore-Penrose pseudoinverse solution with the proposed algorithm in tensors (left), and nested matrices (centre) with its error (right).}
\label{fig:comparison4x2=2}
\end{figure}

Again, as in Figure \ref{fig:comparison2x2=2}, Figure \ref{fig:comparison4x2=2} shows the comparison between the coefficients of the Moore-Penrose pseudoinversa and the ones obtained with the algorithm proposed, in absolute value. The cases shown are three dimensional tensors and nested matrices with again around 100000 data points to validate the results.

\subsubsection{Case 4: $\mathbf{A} : \mathbf{X} =\mathbf{B} \; \rightarrow \; \mathbb{A}_{ijkl}\mathbb{X}_{klmn}=\mathbb{B}_{ijmn}$}

In this case, the order of the unknowns $\mathbb{X}$ is four, and the operation is the double dot product. Since the equation proposed can be automatically adapted to any dimension and both the dot and the double dot product, the algorithm for this case is:\\

\begin{equation}
\widehat{x}_\alpha=\frac{ \mathbb{A}_{\alpha}^T : \left( \mathbb{B}_{\alpha}-\mathbb{A}: \breve{\mathbb{X}}_{\alpha}+\mathbb{A}_{\alpha}\tilde{x}_\alpha\right)}{\mathbb{A}_{\alpha}^T: \mathbb{A}_{\alpha}}
\end{equation}

And in index notation:\\

\begin{equation}
\widehat{\mathbb{X}}_{\hat{\textbf{k}}\hat{\textbf{l}}\hat{\textbf{m}}\hat{\textbf{n}}}=(1-\beta)\tilde{\mathbb{X}}_{\hat{\textbf{k}}\hat{\textbf{l}}}+\beta\frac{ \mathbb{A}_{ij\hat{\textbf{k}}\hat{\textbf{l}}}^T \left(\mathbb{B}_{ij\hat{\textbf{m}}\hat{\textbf{n}}}- \mathbb{A}_{ijop}\breve{\mathbb{X}}_{op\hat{\textbf{m}}\hat{\textbf{n}}}+\mathbb{A}_{ij\hat{\textbf{k}}\hat{\textbf{l}}}\tilde{\mathbb{X}}_{\hat{\textbf{k}}\hat{\textbf{l}}\hat{\textbf{m}}\hat{\textbf{n}}}\right)}{\mathbb{A}_{ij\hat{\textbf{k}}\hat{\textbf{l}}}^T \mathbb{A}_{ij\hat{\textbf{k}}\hat{\textbf{l}}}}
\end{equation}

\begin{figure}[H]
\centering
\subfloat{\includegraphics[width=0.31\textwidth]{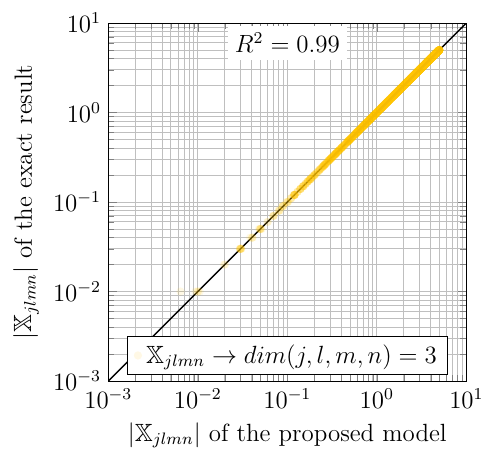}}
\subfloat{\includegraphics[width=0.31\textwidth]{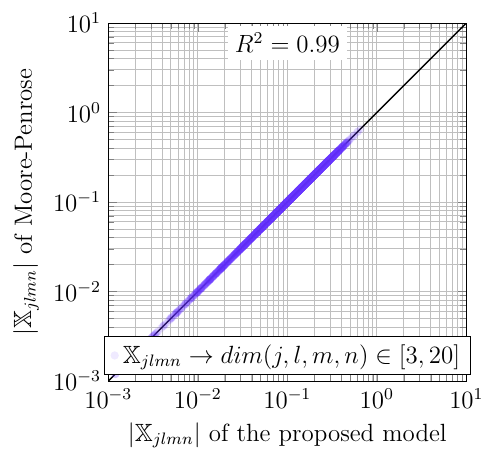}}
\subfloat{\includegraphics[width=0.30\textwidth]{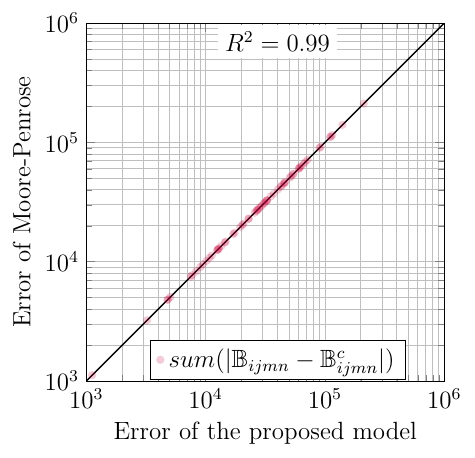}}
\caption{Comparison of the Moore-Penrose pseudoinverse solution with the proposed algorithm in tensors (left), and nested matrices (centre) with its error (right).}
\label{fig:comparison4x4=4}
\end{figure}

As in Figures \ref{fig:comparison2x2=2} and \ref{fig:comparison4x2=2}, Figure \ref{fig:comparison4x4=4} shows the comparison with the Moore-Penrose pseudoinverse in high order tensors and nested matrices. \\

In total, in Figure \ref{fig:comparison4x4=4} there are more than 2 millions of data points with a $R^2$ higher than 0.99, which probes again the agreement between the algorithm proposed and the Moore-Penrose pseudoinverse.\\

\section{Conclusions}

The algorithm proposed has proven to be a generalization of the Gauss-Seidel equation, as well as Jacobi´s, which gives more convergence capabilities due to the consideration of all the terms, not only the diagonal. \\

The exact solution given by Cramer for small systems can also be directly deducted by the algorithm proposed, which proves the validity of its theoretical conceptualization.\\

For larger random systems, where the traditional approach of the Gauss-Seidel equation can not reach a convergence, the algorithm proposed shows a performance comparable to different quadratic solvers. In this line, with around 4 millions of validation points, the algorithm shows its capabilities to reproduce the solution given by the Moore-Penrose pseudoinverse in a wide variety of the random non-square linear systems and under the dot and double dot product operations.

%

\bibliographystyle{plain}
\bibliography{Articulo_algorithm_v9}

\newpage
\section*{Annex: Iterative implementation in Python of the equation \ref{ecuxsysliniter} using the numpy library}
\begin{python}

import numpy as np

# 'matrix of coefficients'
kmat=np.array([......])
# initial guess at iteration 0, 
# and will be the solution afterwards
u=np.array([......])
# vector of constants
b=np.array([......])

convergence=False
# beta for the convergence
beta= ...
# initial error
error0=0.0 
# value of the tolerance for the convergence
tolerance= ... 

while not convergence: # this loop is an iteration

   for j in range(0,len(u),1):
      cj=kmat[:,j]
      kj=np.delete(kmat,j,1)
      uj=np.delete(u,j,0)
      pj=np.dot(kj,uj)
      t1=np.dot(pj-b,cj)
      t2=np.dot(cj,cj)
      
      # updated solution at j
      u0=-(t1/t2)
      up=u[j]
      u[j]=(1-beta)*up+beta*u0
      
   # check convergence
   error=np.linalg.norm(np.dot(kmat,u)-b) 
   if abs(error-error0)<tolerance:
      convergence=True
   else:
      error0=error
      
      
\end{python}

\end{document}